\documentclass[11pt]{article}

\usepackage{amssymb,latexsym,amsmath}

\usepackage{graphicx}

\begin{document}

\newcommand{\bfi}{\bfseries\itshape}

\makeatletter

\@addtoreset{figure}{section}

\def\thefigure{\thesection.\@arabic\c@figure}

\def\fps@figure{h, t}

\@addtoreset{table}{bsection}

\def\thetable{\thesection.\@arabic\c@table}

\def\fps@table{h, t}

\@addtoreset{equation}{section}

\def\theequation{\thesubsection.\arabic{equation}}

\makeatother

\newtheorem{thm}{Theorem}[section]

\newtheorem{prop}[thm]{Proposition}

\newtheorem{lema}[thm]{Lemma}

\newtheorem{cor}[thm]{Corollary}

\newtheorem{defi}[thm]{Definition}

\newtheorem{rk}[thm]{Remark}

\newtheorem{exempl}{Example}[section]

\newenvironment{exemplu}{\begin{exempl}  \em}{\hfill $\surd$

\end{exempl}}

\newcommand{\comment}[1]{\par\noindent{\raggedright\texttt{#1}

\par\marginpar{\textsc{Comment}}}}

\newcommand{\todo}[1]{\vspace{5 mm}\par \noindent \marginpar{\textsc{ToDo}}\framebox{\begin{minipage}[c]{0.95 \textwidth}

\tt #1 \end{minipage}}\vspace{5 mm}\par}

\newcommand{\ea}{\mbox{{\bf a}}}
\newcommand{\eu}{\mbox{{\bf u}}}
\newcommand{\ueu}{\underline{\eu}}
\newcommand{\ueo}{\overline{u}}
\newcommand{\oeu}{\overline{\eu}}
\newcommand{\ew}{\mbox{{\bf w}}}
\newcommand{\ef}{\mbox{{\bf f}}}
\newcommand{\eF}{\mbox{{\bf F}}}
\newcommand{\eC}{\mbox{{\bf C}}}
\newcommand{\en}{\mbox{{\bf n}}}
\newcommand{\eT}{\mbox{{\bf T}}}
\newcommand{\eL}{\mbox{{\bf L}}}
\newcommand{\eV}{\mbox{{\bf V}}}
\newcommand{\eU}{\mbox{{\bf U}}}
\newcommand{\ev}{\mbox{{\bf v}}}
\newcommand{\eve}{\mbox{{\bf e}}}
\newcommand{\uev}{\underline{\ev}}
\newcommand{\eY}{\mbox{{\bf Y}}}
\newcommand{\eK}{\mbox{{\bf K}}}
\newcommand{\eP}{\mbox{{\bf P}}}
\newcommand{\eS}{\mbox{{\bf S}}}
\newcommand{\eJ}{\mbox{{\bf J}}}
\newcommand{\eB}{\mbox{{\bf B}}}
\newcommand{\leb}{{\cal L}^{n}}
\newcommand{\eI}{{\cal I}}
\newcommand{\eE}{{\cal E}}
\newcommand{\hen}{{\cal H}^{n-1}}
\newcommand{\eBV}{\mbox{{\bf BV}}}
\newcommand{\eA}{\mbox{{\bf A}}}
\newcommand{\eSBV}{\mbox{{\bf SBV}}}
\newcommand{\eBD}{\mbox{{\bf BD}}}
\newcommand{\eSBD}{\mbox{{\bf SBD}}}
\newcommand{\ecs}{\mbox{{\bf X}}}
\newcommand{\eg}{\mbox{{\bf g}}}
\newcommand{\paromega}{\partial \Omega}
\newcommand{\gau}{\Gamma_{u}}
\newcommand{\gaf}{\Gamma_{f}}
\newcommand{\sig}{{\bf \sigma}}
\newcommand{\gac}{\Gamma_{\mbox{{\bf c}}}}
\newcommand{\deu}{\dot{\eu}}
\newcommand{\dueu}{\underline{\deu}}
\newcommand{\dev}{\dot{\ev}}
\newcommand{\duev}{\underline{\dev}}
\newcommand{\weak}{\rightharpoonup}
\newcommand{\weakdown}{\rightharpoondown}
\renewcommand{\contentsname}{ }

\title{Microfractured media with a scale and Mumford-Shah energies}
\date{03.2007}
\author{Marius Buliga\footnote{"Simion Stoilow" Institute of Mathematics of the Romanian Academy,
 PO BOX 1-764,014700 Bucharest, Romania, e-mail: Marius.Buliga@imar.ro }}


\maketitle

\begin{abstract}
We want to understand he concentration of damage in microfractured 
elastic media. Due to the different scallings of the volume and area (or area 
and length in two dimensions) the traditional method of homogenization 
using periodic arrays of cells seems to fail when applied to the Mumford-Shah functional 
and to periodically fractured domains. 

In the present paper we are departing from traditional homogenization. The main
result implies the use of Mumford-Shah energies and leads to an explanation of 
the observed  concentration of damage in microfractured elastic bodies. 
\end{abstract}

\section{Introduction}

A new direction of research in brittle fracture mechanics  begins with the article of 
 Mumford \& Shah  \cite{MS} regarding the problem 
of image segmentation. This problem,  which consists in finding the set of edges of a picture 
and constructing 
a smoothed version of that  picture, it turns to be intimately related to the problem of brittle 
crack evolution. 
In the before mentioned article Mumford and Shah propose the
following variational approach to the problem of image segmentation: 
let $g:\Omega \subset\mathbb{R}^{2} \rightarrow [0,1]$ 
be the original picture, given as a distribution of grey levels (1 is white and 0 is black), 
let $u: \Omega \rightarrow R$ 
be the smoothed picture and $K$ be the set of edges. $K$ represents the set where $u$ has jumps, 
i.e. $u \in C^{1}(\Omega \setminus K,R)$. The pair formed by the smoothed picture $u$ and the set 
of edges $K$ minimizes then the functional: 
$$I(u,K) \ = \ \int_{\Omega} \alpha \ \mid \nabla u \mid^{2} \mbox{ d}x \ + \ 
\int_{\Omega} \beta \  \mid u-g \mid^{2} \mbox{ d}x \ + \ \gamma \mathcal{ H}^{1}(K) \ \ .$$
The parameter $\alpha$ controls the smoothness of the new picture $u$, $\beta$ controls the 
$L^{2}$ distance between the smoothed picture and the original one and $\gamma$ controls the 
total length of the edges given by this variational method. The authors remark that for $\beta 
= 0$ the functional $I$ might be useful for an energetic treatment of fracture mechanics. 

An energetic approach to fracture mechanics is naturally suited to explain brittle crack 
appearance under imposed boundary displacements. 
The idea is presented  in the followings. 

The state of a brittle body is  described by a pair 
displacement-crack. $(\eu,K)$ is such a pair if $K$ is a crack --- seen as a surface 
---  
which appears in the body and  $\eu$ is a displacement of the broken body under the imposed 
boundary displacement, i.e. $\eu$ is continuous in the exterior of the surface $K$ and $\eu$
equals the imposed displacement $\eu_{0}$ on the exterior boundary of the body. 

Let us suppose 
that the total energy of the body is a Mumford-Shah functional of the form: 
$$E(\eu,K) \ = \ \int_{\Omega} w(\nabla \eu) \mbox{ d} x \ + \ F(\eu_{0},K) \ \ .$$
The first term of the functional $E$ represents the elastic energy of the body with the
displacement $\eu$. The second term represents the energy consumed to produce the crack 
$K$ in the body, with the boundary displacement $\eu_{0}$   as parameter. 
Then the crack that appears is supposed to be the second term of the pair 
$(\eu,K)$ which minimizes the total energy $E$. 

Models for brittle damage, based on functionals of the Mumford-Shah type have 
have been proposed by Francfort-Marigo \cite{FMa}, Buliga \cite{bu3}, among others. 
Such models have been studied intensively from the mathematical point of view, 
espacially by the Italian school of geometric measure theory, to name a few: 
De Giorgi, Ambrosio, Dal Maso, Buttazzo. 

The first homogenization result, concerning the Mumford-Shah functional, seems 
to be Braides, Defranceschi, Vitali \cite{bradevi}. In this paper it is done 
the homogenization of a Mumford-Shah functional of the form: 
$$\int_{\Omega} f\left(\frac{x}{\varepsilon}, \nabla u\right) \mbox{ d} + 
\int_{S_{u}} g\left( \frac{x}{\varepsilon}, (u^{+}-u^{-} \otimes \nu_{u} \right) 
\mbox{ d}\mathcal{H}^{n-1} \quad .$$ 
The paper Focardi, Gelli \cite{foge} (and the references therein) 
are  part of another 
 line of research which might be relevant for this paper: homogenization of perforated domains. 

In the present paper we are departing from traditional homogenization. The 
line of research concerning perforated domains is close to our problem, but 
for various reasons the results from perforated domains don't apply here. 

We want to understand he concentration of damage in microfractured 
elastic media. Due to the different scallings of the volume and area (or area 
and length in two dimensions) the traditional method of homogenization 
using periodic arrays of cells seems to fail when applied to the Mumford-Shah functional 
and to periodically fractured domains. 

The main result, theorem \ref{thdamage}, implies the use of Mumford-Shah 
energies and leads to an explanation of the observed  concentration of damage in 
microfractured elastic bodies. 

Instead of performing a homogenization 
of the total energy of the microfractured body and then study the minimizers 
of the homogenized energy, we proceed along a different path. We study 
sequences of problems on fractured elastic bodies, indexed by a scale parameter 
$\varepsilon$. Each such problem has (at least approximative) solutions. We find 
estimates of the area of the damaged region in terms of the scale 
$\varepsilon$.

\section{Notations}

Let $\Omega$ be a bounded, open subset of $\mathbb{R}^{2}$, with locally 
Lipschitz boundary. We denote by $\displaystyle Y = [0,1]^{2}$ the unit closed 
square in $\mathbb{R}^{2}$. 

For a given $\varepsilon > 0$ let $\mathbb{Z}_{\varepsilon} \subset \mathbb{R}^{2}$ be the 
lattice of points in $\mathbb{R}^{2}$ with coordinates of the form 
$(\varepsilon m, \varepsilon n)$, for all $m,n \in \mathbb{Z}$. 

We denote by $\displaystyle \mathbb{Z}(\varepsilon, \Omega) \subset \mathbb{Z}_{\varepsilon}$ 
the set of all $\displaystyle z \in \mathbb{Z}_{\varepsilon}$ such that 
$$ z + \varepsilon Y \subset \Omega  \quad . $$
To any $z \in \mathbb{Z}(\varepsilon, \Omega)$ we associate the cell 
$$D_{z} = z + \varepsilon Y \ \subset \ \Omega \quad .$$ 
 
The set  $\displaystyle \mathbb{Z}(\varepsilon, \Omega)$ is finite for any 
$\varepsilon > 0$. We denote the cardinal of this set by $N(\varepsilon)$ and we notice 
that as $\varepsilon$ goes to $0$ we have 
$$\lim_{\varepsilon \rightarrow 0} \frac{N(\varepsilon) \varepsilon^{2}}{A(\Omega)} = 1 \quad , $$ 
where $A(\Omega)$ denotes the area of $\Omega$. Thus for small $\varepsilon$ the 
number of cells $\displaystyle N(\varepsilon)$ is approximately equal to 
$\displaystyle A(\Omega)/ \varepsilon^{2}$. 

\section{The model}

We take $\Omega$ to be the configuration set of a microfractured linear elastic 
body. We explain further what we mean by this. 

The elastic properties of the body are described by an elastic 
potential $$ w: M^{2 \times 2}_{sym}(\mathbb{R}) \rightarrow \mathbb{R} \quad . $$ 
We suppose that the function $w$ is quadratic and strictly positive definite. 

For a given displacement $\displaystyle \eu: \Omega \rightarrow \mathbb{R}^{2}$, the 
elastic energy of the body is given by 
$$\int_{\Omega} w(e(\eu)) \mbox{ d}x \quad , $$
where $e(\eu)$ is the deformation of the displacement $\eu$, that is the symmetric 
part of the gradient of $\eu$: for any $x \in \Omega$ 
$$e(\eu)(x) \ = \ \frac{1}{2} \left( \nabla \eu (x) + \left( \nabla \eu \right)^{T}(x) \right) \quad . $$

For a fixed $\varepsilon > 0$ we suppose that the body contains a distribution 
of micro-fractures at the scale $\varepsilon$, seen as a union of (Lipschitz) 
curves 
$$ F_{\varepsilon} = \bigcup_{z \in \mathbb{Z}(\varepsilon, \Omega)} \left( 
z + \varepsilon F_{z} \right)  \quad , $$
where for each $\displaystyle  z \in \mathbb{Z}(\varepsilon, \Omega)$ the (Lipschitz) 
curve $\displaystyle F_{z}$ lies inside the unit cell $Y$: 
$$F_{z} \subset (0,1)^{2} \quad .$$

We explain further what we mean by an imposed boundary displacement 
$\displaystyle \eu_{0}$, and what we mean by 
 $\eu = \eu_{0}$ on the boundary of $\Omega$. 

We consider, for simplicity, 
that $\eu_{0}: \paromega \rightarrow\mathbb{R}^{n}$ is a continuous and therefore bounded function. Then, 
for any $\eu \in \eSBD(\Omega)$, $\eu = \eu_{0}$  if the approximate limit of $\eu$ equals 
$\eu_{0}$ in any point of $\paromega$  where the first exists, i.e.: 
for all $ x \in \paromega$, if there exists  $ \ev(x)$ such that  
$$\lim_{\rho \rightarrow O_{+}} \frac{\int_{ B_{\rho}(x) \cap \Omega} \mid \eu(y) -  \ev(x) \mid 
\mbox{ d}y }
{\mid B_{\rho}(x) \cap \Omega \mid} = 0 $$
then $\displaystyle \ev(x) \ = \eu_{0}(x)$.

\begin{defi}
The class of admissible  displacements with respect to the distribution of 
cracks $\displaystyle F_{\varepsilon}$ and with respect to the imposed 
displacement $\displaystyle \eu_{0}$ is defined as the  collection of all 
$\displaystyle \eu \in SBD(\Omega)$ such that 
\begin{enumerate} 
\item[(a)] $\displaystyle \eu = \eu_{0}$ on $\paromega$, 
\item[(b)] $\displaystyle F_{\varepsilon} \subset S_{u}$. 
\end{enumerate}
This class of admissible displacements is denoted by 
$\displaystyle Adm(F_{\varepsilon}, \eu_{0})$. 
\label{defadmissible}
\end{defi}

This definition deserves an explanation. An admissible displacement $\eu$ is a 
function which has to be equal to the imposed displacement on the boundary 
of $\Omega$ (condition (a)). Any such function $\eu$ is a special function with 
bounded deformation, that is  a reasonably smooth 
function on the set $\displaystyle \Omega \setminus S_{u}$ and the function $\eu$ 
is allowed to have jumps along the set  $\displaystyle 
S_{u}$. For the technical details see the Appendix. We have to think about 
$\displaystyle S_{u}$ as being a collection of curves, with finite length. 
Physically 
the set $\displaystyle S_{u}$ represents the collection of all cracks in the 
body under the displacement $\eu$. The condition (b) tells us that the 
collection of all cracks associated to  an admissible displacement $\eu$  
contains $\displaystyle F_{\varepsilon}$, at least. 

\begin{defi}
With the notations from definition \ref{defadmissible}, the total energy of 
an admissible displacement $\displaystyle \eu \in Adm(F_{\varepsilon}, \eu_{0})$ 
is given by
$$E_{\varepsilon}(\eu) = \int_{\Omega} w(e(\eu)) \mbox{ d}x \ + \ G \mathcal{H}^{1}(S_{u} \setminus F_{\varepsilon}) \quad . $$
\label{defenergy}
\end{defi} 

The energy of an admissible displacement is of Mumford-Shah type. It contains 
two terms. 

The first term measures the elastic energy of the body under the displacement 
$u$. Notice that in the expression of the elastic energy we have integrated 
over the whole domain $\Omega$. This is simply because the collection of 
cracks associated to $u$ (that is the set $\displaystyle S_{u}$) has Lebesque 
measure $0$, therefore we have 
$$\int_{\Omega} w(e(\eu)) \mbox{ d}x \ = \ 
\int_{\Omega\setminus S_{u}} w(e(\eu)) \mbox{ d}x \quad . $$ 
In physical terms, the right hand side expression would make more sense than 
the left hand side, but from the mathematical point of view they are 
the same. This is not meaning that the elastic energy neglects the fractures. 
Indeed, further we shall infimize the energy $\displaystyle E_{\varepsilon}$ 
over the whole set of admissible displacements. According to condition (b) 
of definition \ref{defadmissible}, this set is defined with respect to the 
collection of cracks $\displaystyle F_{\varepsilon}$, therefore the infimum 
of the energy $\displaystyle E_{\varepsilon}$ depends on the set 
of cracks $\displaystyle F_{\varepsilon}$.

The second term of the Mumford-Shah energy measures the surface energy 
caused by the apparition of new cracks. The collection of new cracks 
is the set $\displaystyle S_{u} \setminus F_{\varepsilon}$. The constant 
$G$ has the dimension of energy per unit area, and it is physically related 
to the Griffith constant.

In  \cite{BCDM} has been proven that functionals like $\displaystyle E_{\varepsilon}$
 are $L^{1}$ inferior semi-continuous and coercive, hence on closed subspaces 
${\eV}$ of $\eSBD(\Omega)$ the functional $\displaystyle E_{\varepsilon}$ 
 has a minimizer.  Such a closed subspace of $\eSBD(\Omega)$ is the space of 
all admissible displacements  $\displaystyle Adm(F_{\varepsilon}, \eu_{0})$. Therefore 
we have: 

\begin{thm}
On the space $\displaystyle Adm(F_{\varepsilon}, \eu_{0})$ we  consider the  
topology given by the convergence: $ \displaystyle \eu_{h}  \rightarrow \eu$  
if
$$ \left\{ \begin{array}{l}
\eu_{h} \ \  L^{2} \ \rightarrow \ \ \eu \ \ , \\
\hen(S_{u_{h}} \Delta S_{u}) \rightarrow 0 \ \ . 
\end{array} \right. $$ 
Then there exists a minimizer of the functional 
$\displaystyle E_{\varepsilon}$ over the set 
$\displaystyle Adm(F_{\varepsilon}, \eu_{0})$. 
\end{thm}

In the following section we shall use approximate minimizers. 

\begin{defi}
For a given $\delta > 0$, a function $\displaystyle \eu \in  
Adm(F_{\varepsilon}, \eu_{0})$ is a $\delta$-approximate minimizer if 
$$E_{\varepsilon}(\eu) \ \leq \ \delta + \inf \left\{ E_{\varepsilon}(\ev) \mbox{ : } 
\ev \in  Adm(F_{\varepsilon}, \eu_{0})\right\} \quad . $$
\label{defaprox}
\end{defi}

For fixed $\delta > 0$, we model an approximate displacement of a microfractured 
body as a sequence of displacements  $\displaystyle \eu_{\varepsilon}$, with 
$\varepsilon$ converging to $0$, such that for each $\varepsilon > 0$ the displacement 
$\displaystyle \eu_{\varepsilon} \in  Adm(F_{\varepsilon}, \eu_{0})$ is a 
$\delta$-approximate minimizer of the Mumford-Shah energy $\displaystyle 
E_{\varepsilon}$, over the set $\displaystyle Adm(F_{\varepsilon}, \eu_{0})$.  

Notice that in the model, at this stage, there is no relation between 
the crack sets $\displaystyle F_{\varepsilon}, F_{\varepsilon'}$, for 
two different scales $\varepsilon, \varepsilon'$. 

\section{An estimate related to damage concentration}

For fixed $\varepsilon, \delta > 0$, given $\displaystyle F_{\varepsilon}$ and 
imposed boundary displacement $\displaystyle \eu_{0}$, let $\displaystyle 
\eu \in  Adm(F_{\varepsilon}, \eu_{0})$ be a $\delta$-approximate 
minimizer of the Mumford-Shah energy $\displaystyle E_{\varepsilon}$. 

In this section we want to estimate the number of $\varepsilon$-cells 
$\displaystyle z + \varepsilon Y$, $z \in \mathbb{Z}(\varepsilon, \Omega)$, 
where the initial cracks $\displaystyle z + \varepsilon F_{z}$ propagated.

Let $l>0$ be a given length. 

\begin{defi}
For any  cell $\displaystyle D_{z} = z + \varepsilon Y$, 
$z \in \mathbb{Z}(\varepsilon, \Omega)$, and any  $\delta$-approximate 
minimizer $\eu$ we define the emergent crack in the cell $\displaystyle D_{z}$ 
by 
  $$S_{u}(z) = \left( z + \varepsilon Y \right) \cap \left( S_{u} \setminus 
\left(z + \varepsilon F_{z}\right) \right) \quad .$$
A cell $\displaystyle D_{z}$ is called active if the length of the emergent 
crack is greater than $\varepsilon l$, that is: 
$$\mathcal{H}^{1}(S_{u}(z)) \geq \varepsilon l \quad . $$
We denote by $M(\varepsilon, l)$ the number of active cells. (In this 
notation we don't mention the dependence of $M(\varepsilon, l)$ on the 
$\delta$-approximate minimizer $\eu$.)
\label{defactive}
\end{defi}

\begin{thm}
Suppose that for fixed $\delta > 0$, the crack sets 
$\displaystyle F_{\varepsilon}$ are chosen so that 
there exists  an approximate displacement of a microfractured 
body  $\displaystyle \eu_{\varepsilon}$, with 
$\varepsilon$ converging to $0$, with the property that the sequence 
$$ \inf \left\{ E_{\varepsilon}(\ev) \mbox{ : } 
\ev \in  Adm(F_{\varepsilon}, \eu_{0})\right\}$$ 
is bounded. 

Then the number of active cells $M(\varepsilon, l)$ is of order 
$1/\varepsilon$ and the area of the damaged region of the body 
$$Damaged(\varepsilon, \Omega) \ = \ \bigcup_{D_{z} \mbox{ active}} D_{z}$$ 
is of order $\varepsilon$.
\label{thdamage}
\end{thm}

\paragraph{Proof.} 
Let $M > 0$ such that for all $\varepsilon > 0$ we have 
$$ \inf \left\{ E_{\varepsilon}(\ev) \mbox{ : } 
\ev \in  Adm(F_{\varepsilon}, \eu_{0})\right\} \leq M \quad .$$ 
According to definition \ref{defaprox}, for any $\varepsilon > 0$ we have 
$$E_{\varepsilon}(\eu_{\varepsilon}) = 
\int_{\Omega} w(e(\eu_{\varepsilon})) \mbox{ d}x + G \mathcal{H}^{1}\left(S_{u_{\varepsilon}} \setminus F_{\varepsilon}\right) \leq $$
$$\leq \delta + \inf \left\{ E_{\varepsilon}(\ev) \mbox{ : } 
\ev \in  Adm(F_{\varepsilon}, \eu_{0})\right\} \leq \delta +  M \quad .$$ 
From definition \ref{defactive} we get the following estimate: 
$$ \mathcal{H}^{1}\left(S_{u_{\varepsilon}} \setminus F_{\varepsilon}\right) =  
\sum_{z \in \mathbb{Z}(\varepsilon, \Omega)}  \mathcal{H}^{1} 
\left( S_{u}(z))\right)  \geq 
 M(\varepsilon, l) \ l \ \varepsilon \quad .$$
We have therefore 
$$ G \ M(\varepsilon, l) \ l \ \varepsilon  \leq G  \mathcal{H}^{1}\left(
S_{u_{\varepsilon}} \setminus F_{\varepsilon}\right) \leq E_{\varepsilon}(\eu_{\varepsilon}) \leq 
M + \delta \quad . $$
All in all we have obtained the estimate: 
$$M(\varepsilon, l) \leq \frac{1}{\varepsilon} \frac{M+\delta}{G l} \quad . $$
The area of the damaged region of the body is 
$$Area(Damaged(\varepsilon, \Omega))  =  \sum_{D_{z} \mbox{ active}} Area(D_{z}) 
 = \varepsilon^{2} M(\varepsilon, l) \leq \varepsilon \frac{M+\delta}{G l} \quad . $$
The proof is done. $\quad \square$

\section{Conclusions}

The theorem implies that the area of the damaged region is much smaller than 
the total area of the body, as $\varepsilon$ goes to zero. In this model the use of Mumford-Shah energies 
leads to an explanation of the observed  concentration of damage in microfractured elastic bodies.

Notice that we need 
more precise estimates in order to prove that the damaged region (at the scale 
$\varepsilon$) converges, as $\varepsilon$ goes to zero, to a curve with 
finite length. All we know at this moment is that the area of the damaged 
region goes to zero as the scale parameter $\varepsilon$. 

In experiments it has been observed that the damaged region is approximately 
straight. It is possible that Mumford-Shah energies might explain this, since 
geometries of the active crack set, that is 
$\displaystyle S_{u}\setminus F_{\varepsilon}$, close to a straight line would 
be preferred by  the energy $\displaystyle E_{\varepsilon}$. See \cite{buligaper} 
for examples that in some situations the leading term of a Mumford-Shah energy 
is the one accounting for the length of the crack, and not the elastic 
energy part. 

Finally, in theorem \ref{thdamage} we obtained an estimate of the number of 
cells where cracks of length at least $\varepsilon \ l$ appear. It would be interested 
to study the interplay between $\varepsilon$ and $l$ in this estimate.

\section{Appendix.Functions with bounded variation or deformation}
\indent

This section is dedicated to a brief voyage trough the spaces $\eSBV$ and 
$\eSBD$. 

The space $\eSBV(\Omega,R^{n})$ of special functions with bounded variation was introduced by 
De Giorgi and Ambrosio in the study of a class of free discontinuity problems 
(\cite{DGA}, \cite{A1}, \cite{A2}). 
For any function $\eu \in L^{1}(\Omega,R^{n})$ let us denote by $D\eu$ the distributional 
derivative of $\eu$ seen as a vector measure. The variation of $D\eu$ is a scalar measure 
defined like this:  
for any Borel measurable subset $B$ of $\Omega$ the variation of $D\eu$ over $B$ is 
\begin{displaymath}
\mid D\eu \mid  (B) \ = \  sup \ \left\{ \sum^{\infty}_{i=1} \mid D\eu  (A_{i}) \mid 
\mbox{ : } \cup_{i=1}^{\infty} A_{i} \subset B \ , \ A_{i} \cap A_{j} = \emptyset \ \ \forall 
i \not = j \right\} \ \ .
\end{displaymath}
A function $\eu$ has bounded variation if the total variation of $D\eu$ is finite. We send the
reader to the book of Evans \& Gariepy \cite{EG} for basic properties of such functions.

 The space $\eSBV(\Omega,R^{n})$ is defined as follows: 
$$\eSBV(\Omega,R^{n}) \ = \ \left\{ \eu \in L^{1}(\Omega,R^{n}) \mbox{ : } \mid D\eu \mid (\Omega)
 < + \infty 
\ , \ \mid D^{s}\eu \mid (\Omega \setminus \eS_{\eu}) = 0 \right\} \ .$$
The Lebesgue set of $\eu$ is the set of points where $\eu$ has approximate limit. The complementary 
set is a $\leb$ negligible set denoted by $\eS_{\eu}$. If $\eu$ is a special function with 
bounded variation then $\eS_{\eu}$ is also $\sigma$ (i.e. countably) rectifiable.  

From the Calderon \& Zygmund \cite{CZ} decomposition theorem we obtain  the following expression of 
$D\eu$, the distributional derivative of $\eu \in \eSBV(\Omega,R^{n})$, seen as a measure: 
$$D\eu \ = \  \nabla \eu (x) \mbox{ d}x \ + \  [\eu ] \otimes \en \mbox{ d}\hen_{|_{K}} \ \ \ .$$ 

We shall use further the notation $\mu \ll \lambda$ if the measure $\mu$ is absolutely continuous 
with respect to the measure $\lambda$. 

Let us define the following Sobolev space  associated to the crack set $K$ (see \cite{ABF}):
$$W^{1,2}_{K} \ = \ \left\{ \eu \in \eSBV(\Omega,R^{n}) \mbox{ : } \int_{\Omega} 
\mid \nabla \eu \mid^{2} 
\mbox{ d}x + \int_{K} [\eu]^{2} \mbox{ d} \hen  <  + \infty \ , \ \mid D^{s}\eu \mid \ll \hen_{|_{K}} 
\right\} \ .$$
It has been proved in \cite{DGCL} the following equality:
\begin{equation}
W^{1,2}(\Omega \setminus K,\mathbb{R}^{n}) 
\cap L^{\infty}(\Omega,R^{n}) \ = \ W^{1,2}_{K}(\Omega,R^{n}) \cap L^{\infty}(\Omega,R^{n}) \ \ .
\label{dg}
\end{equation}

A similar description can be made for the space of special functions with bounded
deformation $\eSBD(\Omega)$ can be found in  \cite{BCDM}. 
For any function $\eu \in L^{1}(\Omega,R^{n})$ we denote by $E\eu$ the symmetric part of 
the distributional derivative of $\eu$, seen as a vector measure. We denote also by $\eJ_{\eu}$ the 
subset of $\Omega$ where $\eu$ has different approximate limits with respect to a
point-dependent direction.   The difference between $\eS_{\eu}$ and $\eJ_{\eu}$ is subtle. Let us 
quote only the fact that for a function $\eu \in \eSBV(\Omega,R^{n})$ the difference of these sets 
is $\hen$-negligible.

The definition of $\eSBD(\Omega)$ 
is the following: 
$$\eSBD(\Omega,R^{n}) \ = \ \left\{ \eu \in L^{1}(\Omega,R^{n}) \mbox{ : } \mid E\eu \mid (\Omega)
 < + \infty 
\ , \ \mid E^{s}\eu \mid (\Omega \setminus \eJ_{\eu}) = 0 \right\} \ .$$
If $\eu$ is a special function with bounded deformation then $\eJ_{\eu}$ is countably rectifiable. 
We have a decomposition theorem for $\eSBD$ functions, similar to Calderon \& Zygmund result
applied  for 
$\eSBV$ functions. The decomposition theorem is due to Belletini, Coscia \& Dal Maso \cite{BCDM}  
and asserts that 
$$E\eu \ = \  \epsilon (\eu) (x) \mbox{ d}x \ + \  [\eu ] \odot \en \mbox{ d}\hen_{|_{\eJ_{\eu}}} \ \ \ .$$ 
Here $\odot$ means the symmetric part of tensor product and $\epsilon (\eu)$ is the approximate 
symmetric gradient, hence the approximate limit of the symmetric part of the gradient of $\eu$. 

We sum up the main facts about functions with bounded variation or deformation, in the 
following three theorems.

\begin{thm}
Let $\eu \in L^{1}(\Omega,\mathbb{R}^{m})$. Then 
\begin{itemize}
\item[-](De Giorgi) If $\eu \in \eBV(\Omega,\mathbb{R}^{m})$ then $\eS_{\eu}$ is countably rectifiable, $\hen(\eS_{\eu} \setminus \eJ_{\eu}) = 0$ and in $\hen$-almost every point  $x \in \eS_{\eu}$
exists the approximate limits of $\eu$ in the directions $\nu(x)$ and $-\nu(x)$ where $\nu(x)$ is the normal to $\eS_{\eu}$ in $x$. 
\item[-](Kohn, Ambrosio, Coscia, Dal Maso) Let $m=n$ and $\eu \in \eBD(\Omega)$. Let $\Theta_{\eu}$ be the Kohn set :
$$\Theta_{\eu} = \left\{ x \in \Omega \mbox{ : } \limsup_{\rho \rightarrow 0^{+}} \frac{ \mid E\eu \mid (B_{\rho}(x))}{\rho^{n-1}} \ > \ 0 \right\}$$
Then $\Theta_{\eu}$ is countably rectifiable , $\eJ_{\eu} \subseteq \Theta_{\eu}$ and $\hen(\Theta_{\eu} \setminus \eJ_{\eu}) = 0$~.
\end{itemize} 
\end{thm}

\begin{thm}
 Let $\eu \in L^{1}(\Omega,\mathbb{R}^{m})$. Then
\begin{itemize}
\item[-](Calderon, Zygmund) If \ $\eu \in \eBV(\Omega,\mathbb{R}^{m})$ \ then \ $\eu$ \ is approximately \ differentiable \newline $\leb$-a.e. in $\Omega$. 
The approximate differential map $x \mapsto \nabla \eu (x)$ is integrable. 
$D\eu$ splits into three mutually singular measures on $\Omega$
$$D\eu = \nabla \eu \mbox{ d}x \ + \ [\eu] \otimes \nu \hen_{|_{\eS_{\eu}}} \ + \ C\eu$$
where $[\eu]$ is the jump of $\eu$ in respect with the normal direction on $\eS_{\eu}$ $\nu$. $C\eu$ is the Cantor part of $D\eu$ defined by $C\eu(A) = D^{s}\eu(A \setminus \eS_{\eu})$ where $D^{s}\eu$ is 
the singular part of $D\eu$ in respect to $\leb$.
\item[-](Belletini, Coscia, Dal Maso)  Let $m=n$ and $\eu \in \eBD(\Omega)$. Then $\eu$ has symmetric approximate differential $\epsilon(\eu)$ $\leb$-a.e. in $\Omega$ and $E\eu$ splits into three mutually
singular measures on $\Omega$
$$E\eu = \epsilon (\eu) \mbox{ d}x \ + \ [\eu] \odot \nu \hen_{|_{\eJ_{\eu}}} \ + \ E^{c}\eu$$ 
Moreover $\eu$ is approximately differentiable $\leb$-a.e. in $\Omega$.
\end{itemize}
\end{thm}

\begin{thm} 
The following are true:
\begin{itemize}
\item[-] $W^{1,1}(\Omega,\mathbb{R}^{m}) \subset \eBV(\Omega,\mathbb{R}^{m})$. The inclusion is continuous in respect with the Banach space topologies. If $$\eu \in \eSBV(\Omega,\mathbb{R}^{m})$$ then 
$$\eu \in W^{1,1}(\Omega \setminus \eS_{\eu},\mathbb{R}^{m})$$ Moreover if $\eu \in W^{1,1}(\Omega \setminus K,\mathbb{R}^{m})\cap L^{\infty}(\Omega , \mathbb{R}^{m})$ , where $K$ is a closed , 
countably rectifiable set with $\hen(K) < + \infty$, then $\eu \in \eSBV(\Omega,\mathbb{R}^{m})$ and $\hen(K \setminus \eS_{\eu}) = 0$.
\item[-] Let $LE^{1}(\Omega)$ be the Banach space of $L^{1}(\Omega,\mathbb{R}^{n})$ functions with $L^{1}$ symmetric differential.  If $\eu \in \eSBD(\Omega)$ then $\eu \in LE^{1}(\Omega\setminus \eJ_{\eu})$.
Let $K$ be a closed , countably rectifiable set with $\hen(K) < + \infty$. If $\eu \in LE^{1}(\Omega\setminus K) \cap  L^{\infty}(\Omega , \mathbb{R}^{n})$ then $\eu \in \eSBD(\Omega)$ and 
 $\hen(K \setminus \eJ_{\eu}) = 0$.
\end{itemize}  
\end{thm}

\end{document}